\documentclass[a4paper,10pt]{article}

\usepackage{amssymb,amsmath,graphics}

\newenvironment{dwd}{\par\noindent{\bf Proof.}}{\par\rightline{$\blacksquare$}}

\newtheorem{theo}{Theorem}
\newtheorem{prop}{Proposition}

\newtheorem{lema}{Lemma}
\newtheorem{defi}{Definition}
\def\be#1\ee{\begin{equation}#1\end{equation}}
\newcommand{\ba}{\begin{eqnarray} }
\newcommand{\ea}{\end{eqnarray} }
\def\bt#1\et{\begin{theo}#1\end{theo}}
\def\bl#1\el{\begin{lema}#1\end{lema}}
\def\bp#1\ep{\begin{prop}#1\end{prop}}
\def\bd#1\ed{\begin{defi}#1\end{defi}}
\global\long\def\sbr#1{\left[#1\right]}
\global\long\def\cbr#1{\left\{#1\right\}}
\global\long\def\rbr#1{\left(#1\right)}
\global\long\def\tbr#1{\left\langle#1\right\rangle}
\global\long\def\abs#1{\left|#1\right|}
\global\long\def\norm#1{\left\lVert#1\right\rVert}
\def\ccA{{\cal A}}

\def\ccN{{\cal N}}

\def\ccP{{\cal P}}

\def\va{\varepsilon}
\def\ra{\rightarrow}

\def\E{\mathbf{E}}

\def\P{\mathbf{P}}
\def\Var{\mathbf{Var}}

\def\N{{\mathbb N}}

\def\R{{\mathbb R}}

\def\ls{\leqslant}
\def\gs{\geqslant}

\def\1{\mathbf{1}}
\def\p{\boldsymbol{p}}

\begin{document}

\title{\bf Some remarks on the Gram-Schmidt walk algorithm and consequences for Komlos conjecture}
\author{{Witold Bednorz and Piotr Godlewski}
\footnote{Research partially supported by NCN Grant UMO-2022/47/B/ST1/02114}
\footnote{Institute of Mathematics, University of Warsaw, Banacha 2, 02-097 Warszawa, Poland}}

\maketitle
\begin{abstract}
In this paper we improve the best known constant for the discrepancy formulated in the Komlos Conjecture. The result is based on the improvement of the subgaussian bound for the random vector constructed in the Gram-Schmidt Random Walk algorithm.
Moreover, we present detailed argument for the smoothed analysis of this random vector.
The analysis concerns a modification of a given matrix in the conjecture by a Gaussian type perturbation. Our result improves the recent paper in this direction.
\end{abstract}

\section{History}

In this paper we discuss some improvement on the constant for the discrepancy problem formulated by Komlos.
Suppose that vectors $v_1,v_2,\ldots,v_n\in \R^d$ and $\norm{v_i}_2\ls 1$ for $i\in[n]$.
The Komlos conjecture concerns the existence of a sequence of signs $(\va_i)^n_{i=1}\in \{-1,1\}^n$
such that
\[
\norm{\sum^n_{i=1}\va_i v_i}_{\infty}\ls C,
\] 
where $C$ is a universal constant. 
First attempts to show the result reaches early eighties, when Beck and Spencer got some partial results towards the conjecture. 
The first important observation was shown in 1985, in papers \cite{Spe1,Spe2} J. Spencer proved the following observation.
\bt
For any $v_1,v_2,\ldots,v_n\in \R^d$ such that $\norm{v_i}_2\ls 1$ there exists a sequence of signs $(\va_i)^n_{i=1}\in \{-1,1\}^n$, such that
\[
 \norm{\sum^n_{i=1} \va_i v_i}_2 \ls C \sqrt{n}\sqrt{1+\log\rbr{\frac{d}{n}}},
\]
where $C$ is a universal constant.
\et
Then, in 1998 W. Banaszczyk \cite{Ban} proved the best up now lower bound for the discrepancy.
\bt
For any $v_1,v_2,\ldots v_n\in \R^d$ there exists
a sequence of signs $(\va_i)^n_{i=1}\in \{-1,1\}^n$
such that
\[
 \norm{\sum^n_{i=1}\va_i v_i}_{\infty}\ls C\sqrt{1+\log d}.
\] 
\et 
Numerous results were proved about the Komlos conjecture to mention the result by Nikolov \cite{Nik} 
who showed some convex-type result that supports the conjecture. The huge success was to discover
the Gram-Schmidt walk algorithm by Bansal et al in 2019 \cite{Bansal}. Quite recently,  in 2022,  there appeared the detail analysis of GSW algorithm by Spielman \cite{Spiel}.

\section{Main results}

In this paper we show how to improve the best known constant for the discrepancy formulated in Komlos conjecture. Our idea is based on the Gram-Schmidt Random Walk (GSW) algorithm which gives a random vector $X=(X_i)^n_{i=1}\in \cbr{\pm 1}^n_{i=1}$
for a certain $d\times n$ matrix $M$ with all columns $M_i\in \R^d$, $1\ls i\ls n$ to be of the Euclidean norm smaller or equal $1$. Variable $X$ is obtained as the result of a certain number of improvements $X_t$, $0\ls t\ls T\ls n$, where $T$ is a random number of running time of the algorithm.  We explain later that some of the GSW algorithm steps may be trivial - they do not give new orthogonal vectors,  that is why we prefer $\hat{T}\ls T$ which represents number of nontrivial time steps. 
The algorithm construction guarantees
\[
\hat{T}\ls \mathrm{dim}\rbr{\mathrm{Lin}\rbr{M_1,M_2,\ldots, M_n}}\ls \min\rbr{d,n}.
\]
Our goal is to replace $d$ in the Banaszczyk lower bound by $\E \hat{T}$.  Namely,
\bt\label{theo1}
For any $v_1,v_2,\ldots,v_n\in \R^d$, there exists a sequence of signs $(X_i)^n_{i=1}\in \{-1,1\}^n$ such that
\[
 \norm{\sum^n_{i=1} v_i X_i}_{\infty}\ls 2\max\rbr{1, \sqrt{2 \E\max_{1\ls i\ls d} Z_{e_i}}\sqrt{\log\E\hat{T}}},
\]
where $Z_{e_i}$ are random variables satisfying $0\ls Z_{e_i}\ls 1$ and $\hat{T}$ is the number of nontrivial times steps of GSW algorithm -  $\hat{T}\ls\min\cbr{d,n}$.
 \et 
Now suppose that $R$ is $d\times n$ matrix of i.i.d. Gaussian entries $R_{i,j}$ of 
$\ccN(0,\frac{\sigma^2}{d})$. Let also $Y$ be a certain modification of $X$ 
distribution, which we detail later, still defined on $\cbr{\pm}^n_{i=1}$. Note that in this direction we follow \cite{Meka}, however our choice of the modification is more natural and closer to the basic definition of $X$. We show that
\bt\label{theo2}
Suppose that $n=\kappa d\log d$, $\kappa,\sigma\gs 1$ and $\va=\sigma \sqrt{\log d}d^{-\kappa/32}$, then
\[
\P_R\rbr{ \P_Y \rbr{\norm{\rbr{M+R}Y}_{\infty}\ls\va}>0}\gs 1-c(d),
\]
where $c(d)\ra 0$, when $d\ra\infty$.
\et 
In fact we can calculate the value of $c(d)$ which is bit complex, what is important however that it behaves like  $1/\log d$.
Our idea can be extended for different $R_{i,j}$ ensambles, however our aim is to give precise lower bound for the basic 
case.

\section{GSW algorithm}

Our improvement is based on the analysis of the GSW algorithm.
Let us quote the algorithm:
\begin{itemize}
 \item Set the initial vector $X_1=\mathbf{0}\in [-1,1]^n$, time step $t\leftarrow 1$, the set of unfrozen indices $\ccA_1\leftarrow [n]$ and the pivot $p_1\leftarrow n$.
 \item Until $\ccA_t\neq \emptyset$ do:
 \begin{enumerate}
 \item Set direction $u_t\leftarrow \mathrm{argmin}_{u\in U}\norm{Mu}_2$, where
 $U=\{u\in\R^n:\;\; u(p_t)=1,\;\; u(i)=0\;\;\forall i\not\in \ccA_t\}$.
 \item Set $\boldsymbol{\delta^{+}}\leftarrow\abs{\max\Delta}$, $\boldsymbol{\delta^{-}}\leftarrow\abs{\min\Delta}$, where
 $\Delta=\{\delta\in\R:\;\; X_t+\delta u_t\in [-1,1]^n\}$.
 \item Let
 \[
 \boldsymbol{ \delta_t }\leftarrow \left\{ \begin{array}{lll}
 \boldsymbol{\delta^{+}} & \mbox{with probab.} & \boldsymbol{\delta^{-}/\rbr{\delta^{-}+\delta^{+}}} \\
 \boldsymbol{\delta^{-}} & \mbox{with probab.} & \boldsymbol{\delta^{+}/\rbr{\delta^{-}+\delta^{+}} }
 \end{array} \right. .
 \]
 \item Set $X_{t+1}\leftarrow X_t+\delta_t u_t$, $\ccA_{t+1}\leftarrow\{i\in [n]:\;\; \abs{X_{t+1}(i)}<1\}$, $p_{t+1}\leftarrow\max{\ccA_{t+1}}$ , $t\leftarrow t+1$.
 \end{enumerate}
 \item Return $X\leftarrow X_{T+1}$, where $\ccA_{T+1}=\emptyset$.
\end{itemize}
Our analysis of the above algorithm starts with an observation that for each time step $\boldsymbol{t}$ we have $\boldsymbol{t}\ls T+1\ls n+1$, where $T$ is random.
We also control pivot number $\p\in [n]$, which indicated the vector $v_{\p}$ used by the algorithm as a pivot.
Note that not all vectors may become the pivot and we denote the set of pivots by $\ccP$. Furthermore, we denote the time step in which $\p$ becomes the pivot by $t_{\p}$ and the pivot number in time step $\boldsymbol{t}$ by $\boldsymbol{p_t}$.
We say that an index $i\in[n]$ is frozen in time step $\boldsymbol{t}$ if $v_{i}\in\ccA_{t}\backslash\ccA_{t+1}$.
\smallskip

\noindent
We define an order $\sigma$ on $[n]$ in the following way. We start by setting $\sigma(n)=\boldsymbol{p_1}=n$. Next, we set $\sigma(n-1)=j$ if $j$ is frozen in time step $\boldsymbol{t}=1$.
If there are more frozen coordinates e.g. $j_1>j_2>\ldots >j_r$ then $\sigma(n-i)=j_i$ for $1\ls i\ls r$. Then we change the time phase and now either we have
a new pivot $\boldsymbol{p_2}\neq \boldsymbol{p_1}$ and then $\sigma(n-r-1)=\boldsymbol{p_2}$
or $\boldsymbol{p_2}=\boldsymbol{p_1}$ and then once again we consider all other frozen coordinates
in this time step and numerate them by $\sigma$ in the same way as described in time step one. 
We continue until time step $T+1$ in which all coordinates are frozen.
\smallskip

\noindent
We now start the orthogonalisation, namely we define orthogonal vectors $w_{\sigma(1)},w_{\sigma(2)},\ldots w_{\sigma(n)}\in\R^d$ by
\[
w_{\sigma(1)}=\frac{v_{\sigma(1)}}{\norm{v_{\sigma(1)}}_2},\;\; w_{\sigma(r)}=
\frac{v_{\sigma(r)}-A_r v_{\sigma(r)}}{\norm{v_{\sigma(r)}-A_rv_{\sigma(r)}}_2},
\]
where $A_r=\sum_{s<r} w_{\sigma(s)}w^{T}_{\sigma(s)}$.
Note that $A_r$ is the orthogonal projection on space
$\mathrm{Lin}(v_{\sigma(1)},\ldots,v_{\sigma(r-1)})$.
By the construction all $w_{\sigma(r)}$ are orthogonal, we have either $\norm{w_{\sigma(r)}}_2=1$ or $w_{\sigma(r)}=0$, and $M u_t=v_{\boldsymbol{p_t}}-A_t v_{\boldsymbol{p_t}}$.
Moreover, $w_{\sigma(1)},\ldots, w_{\sigma(n)}$ form an orthonormal basis of $V=\mathrm{Lin}(v_1,\ldots,v_n)$, thus
\begin{align*}
Mu_t=\sum^n_{r=1}\tbr{w_{\sigma(r)},v_{\boldsymbol{p_t}}-A_tv_{\boldsymbol{p_t}}}w_{\sigma(r)}=\sum^n_{r=1}\sbr{\tbr{w_{\sigma(r)},v_{\boldsymbol{p_t}}}-\tbr{A_tw_{\sigma(r)},v_{\boldsymbol{p_t}}}}w_{\sigma(r)}.
\end{align*}
Let $\ell_t=\abs{\ccA_t}$. It is crucial to observe that $A_t w_{\sigma(r)}=w_{\sigma(r)}$ if $r< \ell_t$ since $w_{\sigma(r)}$ is already in the space on which $A_t$ projects.
On the other hand, by the construction $A_t w_{\sigma(r)}=0$ for all $r\gs\ell_t$. Finally let $g_t=\sigma^{-1}(\boldsymbol{p_t})$, for any $r>g_t$ we have $w_{\sigma(r)}\bot v_{\boldsymbol{p_t}}$.
Therefore
\begin{equation}\label{eq1.1}
M u_t=\sum^{g_t}_{r=\ell_t}\tbr{w_{\sigma(r)},v_{\boldsymbol{p_t}}}w_{\sigma(r)}.
\end{equation}
\smallskip

\noindent Let us fix pivot $\p$ and let
\[
 I_{\p}=\{t\gs 1:\;\; \boldsymbol{p_t}=\p\}=\{t\gs t_{\p}:\;\; \boldsymbol{p_t}=\p\}
\]
be the set of time steps for which $\p$ is the pivot.
There are $k+1$ vectors frozen while $\p$ is the pivot, i.e.:
\begin{equation}\label{eq1.2}
 w_{\sigma(r-k)}, w_{\sigma(r-k-1)},\ldots, w_{\sigma(r)},
\end{equation}
where
\[
 r=\sigma^{-1}(\p),\;\;k=\abs{\ccA_{t_{\p}}\backslash\ccA_{t_{\p}+\abs{I_{\p}}}}-1\leq n-t_{\p}.
\]
We can split $V$ into orthogonal subspaces $\rbr{V_{\p}}_{\p\in\ccP}$, where $V_{\p}$ is spanned by \eqref{eq1.2}. Furthermore, we can define the main projection matrix $P$ on $V$ and projection matrices $P_{\p}$ with respect to each subspace $V_{\p}$. More precisely
\begin{equation}\label{eq1.3}
 P=\sum_{\p\in\ccP}P_{\p},\;\; P_{\p}=\sum^k_{s=0}w_{\sigma(r-s)}w^{T}_{\sigma(r-s)}.
\end{equation}
Now we split the set of coordinates $r,r-1,\ldots,r-k$ into sets $Q^{\p}_{t_{\p}-1},Q^{\p}_{t_{\p}},\ldots,Q^{\p}_{n}$,
where
\begin{align*}
 Q^{\p}_{t}=
 \begin{cases}
 \cbr{\sigma^{-1}(\p)}=\cbr{r} & \text{if } t=t_{\p}-1 \\
 \cbr{\sigma^{-1}(j):\;\; j\in\ccA_{t}\backslash\ccA_{t+1},\;j\neq\p} & \text{if } t\in I_{\p} \\
 \emptyset & \text{otherwise}
 \end{cases}.
\end{align*}
This allows us to split $V_{\p}$ further into orthogonal subspaces $V_{Q^{\p}_{t}}$, spanned by vectors from $Q^{\p}_{t}$, with corresponding projection matrices
\begin{equation}\label{eq1.4}
    P_{Q^{\p}_{t}}=\sum_{r\in Q^{\p}_t}w_{\sigma(r)}w^{T}_{\sigma(r)}.
\end{equation}
For any $t\gs t_{\p}$ we can now rewrite \eqref{eq1.1} as
\begin{equation}\label{eq1.5}
M u_t=\sum^{t-1}_{s=t_{\p}-1}\sum_{r\in Q^{\p}_s}\tbr{w_{\sigma(r)},v_{\p}}w_{\sigma(r)}.
\end{equation}
We want to estimate the change in discrepancy for pivot $\p$, i.e.
\[
\sum_{t\in I_{\p}}M\rbr{X_{t+1}-X_{t}}=\sum_{t\in I_{\p}}M\delta_t u_t.
\]
Using \eqref{eq1.5} we obtain the following formula
\begin{equation}\label{eq1.6}
\sum_{t\in I_{\p}}\delta_t\tbr{Mu_t,v}=\sum^n_{t=t_{\p}}\delta_t\sum^{t-1}_{s=t_{\p}-1}\sum_{r\in Q^{\p}_s}\tbr{w_{\sigma(r)},v_{\p}}\tbr{w_{\sigma(r)},v}.
\end{equation}

\section{Analysis of the Laplace transform}

Let $\Delta_{\p}$ be generated by $\delta_t$ for $t<t_{\p}$. We define
\[
\alpha_r=\tbr{w_{\sigma(r)},v_{\p}},\;\;\beta_r=\tbr{w_{\sigma(r)},v}.
\]
and
\[
F=\frac{1}{2}\rbr{\sum^{n}_{t=t_{\p}-1}\abs{\sum_{r\in Q^{\p}_t}\alpha_r\beta_r}}^2.
\]
We show that conditionally $F$ is the right
\begin{prop}\label{pro1}
The following inequality holds true
\[
\E\rbr{\exp\rbr{\sum_{t\in I_{\p}}\delta_t\tbr{Mu_t,v}-F} \;\middle|\; \Delta_{\p}}\ls 1.
\]
\end{prop}
\begin{dwd}
For any $R=t_{\p},\ldots,n$ we define
\[
S_{\p}(R)=\sum^n_{t=t_{\p}}\delta_t\sum^{(R\wedge t)-1}_{s=t_{\p}-1}\sum_{r\in Q^{\p}_s}\alpha_r\beta_r.
\]
Moreover, let
\[
F_{\p}(R)=\frac{1}{2}\sbr{\sum^{R-1}_{t=t_{\p}-1}\abs{\sum_{r\in Q^{\p}_t}\alpha_r\beta_r}}^2
\]
and
\[
g(R)=\E\rbr{\exp\rbr{S_{\p}(R)-F_{\p}(R)} \;\middle|\; \Delta_{\p}}.
\]
We are going to show that
\[
g(n+1)\gs g(n)\ls g(n-1)\ls\ldots g(t_{\p})\ls 1.
\]
The inequality $g(n+1)\gs g(n)$ is trivial. 
We aim to show that for $R+1\ls n$, it holds true that
\begin{equation}\label{eq2.1}
g(R+1)=\E\rbr{\exp\rbr{S_{\p}(R+1)-F_{\p}(R+1)} \;\middle|\; \Delta_{\p}}\ls\E\rbr{\exp\rbr{S_{\p}(R)-F_{\p}(R)} \;\middle|\; \Delta_{\p}}=g(R).
\end{equation}
We observe that for $R<n$ we have
\begin{align*}
\E\rbr{\exp\rbr{S_{\p}(R+1)-F_{\p}(R+1)} \;\middle|\; \Delta_{\p}(R)}&=
\exp\rbr{\sum^{R-1}_{t=t_{\p]}}\delta_t\sum^{t-1}_{s=t_{\p}-1}\sum_{r\in Q^{\p}_s}\alpha_r\beta_r-F_{\p}(R+1)}\\
&\cdot\E\exp\rbr{\sum^n_{t=R}\delta_t\sum^R_{s=t_{\p}-1}\sum_{r\in Q^{\p}_s}\alpha_r\beta_r \;\middle|\; \Delta_{\p}(R)},
\end{align*}
where $\Delta_{\p}(R)$ is generated by all $\delta_t$ for $t\ls R-1$.
Clearly, $\exp\rbr{S_{\p}(R)-F_{\p}(R)}$ is $\Delta_{\p}(R)$-measurable.
Therefore, to show \eqref{eq2.1} is suffices to show that
\begin{equation}\label{eq2.2}
\E\exp\rbr{\sum^n_{t=R}\delta_t\sum^R_{s=t_{\p}-1}\sum_{r\in Q^{\p}_s}\alpha_r\beta_r \;\middle|\; \Delta_{\p}(R)}\ls\exp(F_{\p}(R+1)-F_{\p}(R)).
\end{equation}
Note that under $\Delta_{\p}(R)$ condition all $\alpha_r,\beta_r$ for $r\in Q^{\p}_s$, $s\ls R$
are fixed. Moreover, variable $\delta_R=\sum^n_{t=R}\delta_t$ under $\Delta_{\p}(R)$ can attain only two values: $1\pm z_{\p}(R)$, where $\abs{z_{\p}(R)}\ls 1$ and $\E\rbr{\delta_R \;\middle|\; \Delta_{\p}(R)}=0$. More precisely there exists $z_{\p}(R)\in [-1,1]$ which is $\Delta_{\p}(R)$-measurable such that
\begin{equation}\label{eq2.3}
\begin{array}{lll}
\P(\delta_R=-1-z_{\p}(R)) &=&\frac{1-z_{\p}(R)}{2}\\
\P(\delta_R=1-z_{\p}(R)) &=&\frac{1+z_{\p}(R)}{2}
\end{array}.
\end{equation}
Consequently, the proof of \eqref{eq2.2} reduces to the following result.
\begin{lema}\label{lemi1}
For any $x\in [-1,1]$ and $a,b\in \R$ it holds true that
\[
\frac{1-x}{2}\exp(-(1+x)(a+b))+\frac{1+x}{2}\exp((1-x)(a+b))\ls\exp\rbr{\abs{a}\abs{b}+\frac{1}{2}\abs{b}^2}.
\]
\end{lema}
\begin{dwd}
Let us define
\[
L=\frac{(1-x)e^{-(1+x)(a+b)}+(1+x)e^{(1-x)(a+b)} }{(1-x)e^{-(1+x)a}+(1+x)e^{(1-x)a}}.
\]
Observe that
\[
L=pe^{-b(1+x)}+qe^{b(1-x) }, 
\]
where
\[
p=\frac{(1-x)e^{-a}}{(1-x)e^{-a}+(1+x)e^{a}},\;\;q=\frac{(1+x)e^{a}}{(1-x)e^{-a}+(1+x)e^{a}}.
\]
Obviously,
\[
S=-p(1+x)b+q(1-x)b=\rbr{1-x^2}\frac{e^a-e^{-a}}{e^{a}+e^{-a}+x\rbr{e^a-e^{-a}}}b=\frac{(1-x^2)\sinh a}{\cosh{a}+x\sinh{a}}b.
\]
We show that
\[
\abs{\frac{(1-x^2)\sinh a}{\cosh{a}+x\sinh{a}}}\ls F(x)=\frac{(1-x^2)\abs{\sinh{a}}}{\cosh{a}-\abs{x}\abs{\sinh{a}}}\ls\abs{a}.
\]
The maximum of $F$ is attained at
\[
\abs{x}=\frac{\cosh{a}-1}{\abs{\sinh{a}}}
\]
and is equal
\[
\frac{2\rbr{\cosh{a}-1}}{\abs{\sinh{a}}}.
\]
We end up in the well known inequality
\[
\cosh{a}-1\ls\frac{1}{2}\abs{a}\abs{\sinh{a}}
\]
which proves that $\abs{S}\ls\abs{a}\abs{b}$. We can now use Hoeffding lemma to show that
\[
Le^{-S}=pe^{-(1+x)b-S}+qe^{(1-x)b-S}\ls e^{\frac{1}{8}(2b)^2}=e^{\frac{1}{2}b^2}.
\]
Hence
\[
L\ls e^{S+\frac{1}{2}b^2}\ls e^{\abs{a}\abs{b}+\frac{1}{2}\abs{b}^2}.
\]
\end{dwd}
The Bellman approach presented above shows that
$g(n+1)\ls g(t_{\p})$. The last step is to observe that by Hoeffding lemma
\begin{equation}\label{eq2.4}
g(t_{\p})=\E\rbr{\exp\rbr{\sum^n_{t=t_{\p}}\delta_t\sum_{r\in Q^{\p}_{t_{\p}-1}}\alpha_r\beta_r-\frac{1}{2}\abs{\sum_{r\in Q^{\p}_{t_{\p}-1}}\alpha_r\beta_r}^2} \;\middle|\; \Delta_{\p}}\ls 1.
\end{equation}
Note that $Q^{\p}_{t_{\p}-1}=\{\sigma^{-1}(\p)\}$, thus $\sum_{r\in Q^{\p}_{t_{\p}-1}}\alpha_r\beta_r$ is $\Delta_{\p}$-measurable. 
Moreover, similarly as stated in \eqref{eq2.3},
random variable $\delta_{t_{\p}}=\sum^n_{t=t_{\p}}\delta_t$ under $\Delta_{\p}$
condition can attain only two values: $1\pm z_{\p}(t_{\p})$, where $\abs{z_{\p}(t_{\p})}\ls 1$.
We recall that $\E\rbr{\delta_{t_{\p}} \;\middle|\; \Delta_{\p}}=0$. This is due to the construction, there exists $z_{\p}(t_{\p})\in [-1,1]$ which is $\Delta_{\p}$-measurable such that
\[
\begin{array}{lll}
\P(\delta_{t_{\p}}=-1-z_{\p}(t_{\p})) & =& \frac{1-z_{\p}(t_{\p})}{2} \\
\P(\delta_{t_{\p}}=1-z_{\p}(t_{\p})) & = & \frac{1+z_{\p}(t_{\p})}{2}
\end{array}.
\]
Thus \eqref{eq2.4} follows from Hoeffding lemma.
This ends the proof.
\end{dwd}
The last step is to use the above result in the induction procedure.
Let us define
\[
Z_v=\sum_{\p\in\ccP}\rbr{\sum^{n}_{s=t_{\p}-1}\abs{\sum_{r\in Q^{\p}_s}\alpha_r\beta_r}}^2.
\]
Note that $Q^{\p}_s=\emptyset$ if $\boldsymbol{p_s}\neq\p$.
Using the global induction we have proved the following result.
 \bt
The following inequality holds
\[
\E\rbr{\exp\rbr{\tbr{MX,v}-\frac{1}{2}Z_v}}\ls 1.
\]
Moreover $Z_v\ls\norm{v}^2$.
\et
The main improvement is that we may now refer to $Z_{e_i}$ for $i\in [d]$.
Let $L_{\p}$ be the number of nonempty $Q^{\p}_t$ for $t=t_{\p}-1,t_{\p}, \ldots, n$.
Note that $\hat{T}=\sum_{\p\in\ccP}L_{\p}\ls m$, $m=\mathrm{dim}(V)$,
since $Q^{\p}_t$ are pairwise disjoint and $v_r$, $r\in Q^{\p}_t\neq\emptyset$ form a linear basis of $V$.
Hence, using projection matrix formulas \eqref{eq1.3} and \eqref{eq1.4},
\begin{align*}
\sum^{d}_{i=1}Z_{e_i}&=\sum^d_{i=1}\sum_{\p\in\ccP}\rbr{\sum^{n}_{s=t_{\p}-1}\abs{\sum_{r\in Q^{\p}_s}\tbr{v_{\p},w_{\sigma(r)}}\tbr{w_{\sigma(r)},e_i}}}^2\\
&\ls\sum_{\p\in\ccP}\sum^d_{i=1}\rbr{\sum^{n}_{s=t_{\p}-1}\abs{v_{\p}P_{Q^{\p}_s}e_i}}^2\ls\sum_{\p\in\ccP}\sum^d_{i=1}L_{\p}\sum^{n}_{s=t_{\p}-1}\abs{v_{\p}P_{Q^{\p}_s}e_i}^2\\
&\ls\sum_{\p\in\ccP}L_{\p}\norm{P_{\p}v_{\p}}_2^2\ls\sum_{\p\in\ccP}L_{\p}=\hat{T}\ls m.
\end{align*}
Let us observe that
\begin{equation}\label{eq2.5}
\E \rbr{\exp\rbr{\lambda\tbr{MX,v}-\frac{1}{2}\lambda^2 Z_v}}\ls 1.
\end{equation}
Note that $(X_t)_{t\gs 1}$ is a martingale, so
\begin{equation}\label{eq2.6}
\E\tbr{MX,v}=0.
\end{equation}
Moreover,
\[
\E\cosh\rbr{\lambda\tbr{MX,v}}\exp\rbr{-\frac{1}{2}\lambda^2 Z_v}\ls 1
\]
and therefore, by Jensen's and $e^{-x}\geq 1-x$ inequalities, we obtain
\begin{align}\label{eq2.7}
\begin{split}
\E\rbr{\cosh\rbr{\lambda\tbr{MX,v}}-1}\exp\rbr{-\frac{1}{2}\lambda^2 Z_v}&\ls\E\rbr{1-\exp\rbr{-\frac{1}{2}\lambda^2 Z_v}}\\
&\leq 1-\exp\rbr{-\frac{1}{2}\lambda^2\E Z_v}\ls\frac{1}{2}\lambda^2\E Z_v.
\end{split}
\end{align}
Let us fix a constant $C>0$. Note that we can split event $\cbr{\norm{MX}_{\infty}>C}$ into sets $A_1,A_2,\ldots A_m$ which are pairwise 
disjoint and 
\[
A_i=\cbr{\abs{\tbr{MX,e_i}}>C}\cap\bigcap^{i-1}_{j=1}\cbr{\abs{\tbr{MX,e_i}}\ls C}.
\]
Let us observe that
\[
\exp\rbr{-\frac{1}{2}\lambda^2 \E Z_{e_i}\1_{A_i} /\P(A_i)}\ls \E \exp\rbr{-\frac{1}{2} \lambda^2 Z_{e_i}}\1_{A_i} /\P(A_i).
\]
Suppose that $\P\rbr{\norm{MX}_{\infty}>C}=1$, which means $\sum^d_{i=1}\P(A_i)=1$ and by Jensen's inequality and \eqref{eq2.7} we have
\begin{align*}
\exp\rbr{-\frac{1}{2}\lambda^2\sum^d_{i=1}\E Z_{e_i}\1_{A_i}}&\ls\sum^d_{i=1}\P(A_i)\exp\rbr{-\frac{1}{2}\lambda^2\E Z_{e_i}\1_{A_i}/\P(A_i)}\ls\sum^d_{i=1}\E\exp\rbr{-\frac{1}{2}\lambda^2 Z_{e_i}}\1_{A_i}\\
&\ls\frac{1}{\cosh\rbr{C\lambda}-1}\sum^d_{i=1}\E \rbr{\cosh\rbr{\lambda\tbr{MX,e_i}}-1}
\exp\rbr{-\frac{1}{2}\lambda^2 Z_{e_i}}\\
&\ls\frac{\lambda^2}{2\rbr{\cosh\rbr{C\lambda}-1}}\sum^d_{i=1}\E Z_{e_i} \ls \frac{\lambda^2 \E \hat{T}}{2\rbr{\cosh\rbr{C\lambda}-1}}.
\end{align*}
On the other hand,
\[
\sum^d_{i=1}\E Z_{e_i}\1_{A_i}\ls \E \max_{1\ls i\ls d} Z_{e_i}\ls 1.
\]
Therefore 
\[
\exp\rbr{-\frac{1}{2}\lambda^2 \E\max_{1\ls i\ls d }Z_{e_i}}\ls \frac{\lambda^2}{2\rbr{\cosh\rbr{C\lambda}-1}}.
\]
Note that 
\[
\frac{\lambda^2}{2\rbr{\cosh\rbr{C\lambda}-1}}\ls \frac{\lambda^2}{e^{C\lambda}-1-C\lambda}.
\]
Thus for $C\gs 2$ the following inequality holds true
\[
e^{C\lambda}-1 -C\lambda > \lambda^2 e^{C\lambda/2}
\]
since
\begin{align*}
\sum^{\infty}_{k=2}\frac{\rbr{C\lambda}^k}{k!}=\sum^{\infty}_{l=0}\frac{\rbr{C\lambda}^{l+2}}{(l+2)!}=\sum^{\infty}_{l=0}\frac{\rbr{C\lambda/2}^{l}}{l!}\frac{2^l\rbr{C\lambda}^2}{(l+1)(l+2)}>\lambda^2 e^{C\lambda/2}.
\end{align*}
We have used here that $C\gs2$ and $2^{l+2}>(l+1)(l+2)$ for all $l\in\N$. 
Hence
\[
\exp\rbr{-\frac{1}{2}\lambda^2 \E\max_{1\ls i\ls d} Z_{e_i}}<\exp\rbr{-\frac{1}{2}C\lambda}\E\hat{T}.
\]
The optimal $\lambda$ equals $C/\rbr{2\E\max_{1\ls i\ls d }Z_{e_i}}$ and then
\[
\frac{1}{2}C\lambda-\frac{1}{2}\lambda^2 \E\max_{1\ls i\ls d } Z_{e_i} = \frac{C^2}{8\E\max_{1\ls i\ls d} Z_{e_i}}
\]
So if $C\gs 2$ then
\[
\exp\rbr{\frac{C^2}{8\E\max_{1\ls i\ls d} Z_{e_i}}}<\E\hat{T}
\]
so
\[
C < 2\sqrt{2 \E\max_{1\ls i\ls d} Z_{e_i}} \sqrt{\log\E\hat{T}}.
\]
Thus if $C\gs 2\max\cbr{1,\sqrt{2 \E\max_{1\ls i\ls d} Z_{e_i}} \sqrt{\log \E \hat{T}}}$, we have a contradiction. This completes the proof.
\smallskip

\noindent
This shows the main improvement of the Banaszczyk result
\begin{theo}
There exists $X=(X_i)^n_{i=1}\in \{-1, 1\}^n$ such that
\[
\norm{MX}_{\infty}\ls 2\max\rbr{1,\sqrt{2 \E\max_{1\ls i\ls d} Z_{e_i}} \sqrt{\log\E\hat{T}}},
\]
where $K$ is a universal constant.
\end{theo}

\section{The smoothed analysis}

We now turn to the analysis of a perturbation of $M$. Suppose that $\bar{X}$ 
is generated by the GSW algorithm however not for the original matrix $M$ but rather for the matrix $\bar{M}=2^{-1/2}(M^T,\mathrm{Id}_n)^T$, where $\mathrm{Id}_n$ is the identity matrix $n\times n$. By our construction
\[
\E\exp\rbr{\lambda\tbr{\bar{M}\bar{X},(w,v)}-\frac{1}{2}\lambda^2 Z_{(w,v)}}\ls 1.
\]
Thus, we have two properties 
\begin{equation}\label{property1}
\E \exp\rbr{\lambda\tbr{M\bar{X},w}-\lambda^2 Z_w }\ls 1,
\end{equation} 
which in particular means
\[
\E\tbr{M\bar{X},w}^2 \ls 2 \E Z_w,\;\;\mbox{hence}\;\; \E\norm{M\bar{X}}_2^2=2V\ls 2\E T.
\]
Moreover,
\begin{equation}\label{property2}
\E\exp\rbr{\lambda\tbr{\bar{X},v}}\ls\exp\rbr{\lambda^2\norm{v}^2},
\end{equation} 
since $Z_v\ls\norm{v}^2$.
\smallskip

\noindent
The point is that we have to slightly modify the distribution of $\bar{X}$. More precisely, we use
$X$ distributed on $\{\pm 1\}^n$, however for any $A\subset \{\pm 1\}^n$
\[
\P(X\in A)=W^{-1}\E \1_{\bar{X}\in A} \1_{\bar{X}\in B }\exp\rbr{\frac{d\norm{M\bar{X}}_2^2}{2\sigma^2 n }},\;\;W= \E \1_{\bar{X}\in B }
\exp\rbr{\frac{d\norm{M\bar{X}}_2^2}{2\sigma^2 n}},
\]
where $\sigma^2$ is a parameter and
\[
B=\cbr{x\in \cbr{\pm 1}^n:\;\; \norm{Mx}_2^2\ls 2C V}\;\; \mbox{and} \;\; 2V=\E \norm{M\bar{X}}_2^2.
\]
Note that
\begin{equation}\label{property3}
2V=\E\norm{M\bar{X}}_2^2=\sum^d_{i=1}\Var\rbr{\tbr{MX,e}}\ls 2\E\hat{T}\ls 2d,
\end{equation}
however much better upper bound for variance is proved in \cite{Spiel}.
\smallskip

\noindent
Suppose now there is a random $d\times n$ matrix $R$ independent of $X$, where rows $R_i$, $i\in [d]$
are independent. We require that $R_{i,j}$ are independent of Gaussian distribution $\ccN(0,\sigma/\sqrt{d})$.
Note that conditionally on $X$ variables $\tbr{M_i+R_i,X}$ are independent. 
We are going to use the second moment approach, namely
\begin{align*}
& \P_R\rbr{ \P_X \rbr{\norm{\rbr{M+R}X}_{\infty}\ls\va}>0} \\
& \gs \frac{\sbr{{\E_R \P_X \rbr{\norm{\rbr{M+R}X}_{\infty}\ls\va}}}^2}{\E_R \sbr{ \P_X \rbr{\norm{\rbr{M+R}X}_{\infty}\ls\va}}^2}.
\end{align*}
Our main result states that if $n\sim d\log d$, and $\va=d^{-c}$, where $c\sim 1/4$
then the above probability is almost $1$.
Formally, we aim to show Theorem \ref{theo2} with $X$ as the required modification.
\smallskip

\noindent
Consider $X,Y$ independent of the same distribution. 
We observe that
\begin{align*}
& \E_R \sbr{ \P_X\rbr{\cbr{\norm{\rbr{M+R}X}_{\infty}\ls\va}}}^2\\ 
&= \E_{X,Y} \P_R \rbr{\bigcap^d_{i=1} \cbr{\abs{\tbr{M_i+R_i,X}}\ls\va}\cap\bigcap^d_{i=1}\cbr{\abs{\tbr{M_i+R_i,Y}}\ls\va}} \\
& = \E_{X,Y} \prod^d_{i=1} \P_R\rbr{ \cbr{\abs{\tbr{M_i+R_i,X}}\ls\va}\cap\cbr{\abs{\tbr{M_i+R_i,Y}}\ls\va}}.
\end{align*}
On the other hand, 
\begin{align*}
&\sbr{\E_R \P_X\rbr{\norm{\rbr{M+R}X}_{\infty}\ls\va}}^2\\
&=\E_{X,Y} \P_{R^1}\rbr{\bigcap^d_{i=1}\abs{\tbr{M_i+R^1_i,X}}\ls\va}
\P_{R^2}\rbr{\bigcap^d_{i=1}\abs{\tbr{M_i+R^2_i,X}}\ls\va}\\
&=\E_{X,Y} \prod^d_{i=1} \P_{R^1}\rbr{\abs{\tbr{M_i+R^1_i,X}}\ls\va}\P_{R^2}\rbr{\abs{\tbr{M_i+R^2_i,Y}}\ls\va} .
\end{align*}

\section{Main comparison result}

We are going to show that that there exists $C_i(X,Y)$ for each $i\in d$
\begin{align}
 & \P_R\rbr{\{\abs{\tbr{M_i+R_i,X}}\ls\va\}\cap\{\abs{\tbr{M_i+R_i,Y}}\ls\va\}}\nonumber\\
 & \ls 
 C_i(X,Y)\P_{R^1}\rbr{\abs{\tbr{M_i+R^1_i,X}}\ls\va}\P_{R^2}\rbr{\abs{\tbr{M_i+R^2_i,Y}}\ls\va} \label{ineq3.1}.
\end{align}
Let us now fix $X=x\in \cbr{\pm 1}^n_{j=1}$. We choose $R_{ij}=\frac{\sigma}{\sqrt{d}} G_{ij}$, where $G_{ij}$ are independent symmetric Gaussian variables. 
Clearly for each $i\in [d]$ we have
\begin{align*}
& \P_{G_i}\rbr{\abs{\frac{\sigma}{\sqrt{d}}\sum^n_{j=1}G_{i,j}x_j-\sum^n_{j=1}M_{ij}x_j}\ls\va,\;\; \abs{\frac{\sigma}{\sqrt{d}} \sum^n_{j=1}G_{ij}y_j-\sum^n_{j=1}m_{ij}y_j}\ls\va}\\
&=\P_{g_1,g_2}\rbr{\cbr{\abs{\frac{\sigma}{\sqrt{d}}\rbr{\sqrt{k}g_1+\sqrt{n-k}g_2}-m_i(1)}\ls\va}
\cap\cbr{\abs{\frac{\sigma}{\sqrt{d}}\rbr{\sqrt{k}g_1-\sqrt{n-k}g_2}-m_i(2)}\ls\va}}\\
&=\P_{g_1,g_2}\rbr{\abs{\frac{\sigma}{\sqrt{d}}\sqrt{k}g_1-\frac{1}{2}\rbr{m_i(1)+m_i(2)}}+\abs{\frac{\sigma}{\sqrt{d}}\sqrt{n-k}g_2-\frac{1}{2}\rbr{m_i(1)-m_i(2)}}\ls\va},
\end{align*}
where $g_1,g_2$ independent Gaussian r.v.'s and
\[
m_i(1)=\tbr{M_i,x},\;\; m_i(2)=\tbr{M_i,y}.
\]
We denote
\begin{equation}\label{eq3.2}
k=\abs{\cbr{j\in [n]:\;\; x_i=y_i}},\;\;n-k=\abs{\cbr{j\in [n]:\;\; x_i=y_i}}.
\end{equation}
Note that
\[
k=\frac{1}{4}\abs{x+y}^2,\;\;n-k=\frac{1}{4}\abs{x-y}^2, \;\;2k-n=\tbr{x,y}.
\]
We have
\begin{align*}
&\P_{g_1,g_2}\rbr{\abs{\frac{\sigma}{\sqrt{d}}\sqrt{k}g_1-\frac{1}{2}\rbr{m_i(1)+m_i(2)}}+\abs{\frac{\sigma}{\sqrt{d}}\sqrt{n-k}g_2-\frac{1}{2}\rbr{m_i(1)-m_i(2)}}\ls\va}\\
&=\frac{1}{2\pi}\frac{d}{\sigma^2}\frac{1}{\sqrt{k}\sqrt{n-k}}\int_{\abs{\bar{u}}+\abs{\bar{v}}\ls\va}
\exp\rbr{-\frac{d}{2\sigma^2 k}\abs{\bar{u}-\frac{1}{2}\rbr{m_i(1)+m_i(2)}}^2}\\
&\cdot\exp\rbr{\frac{d}{2\sigma^2 (n-k)}\abs{\bar{v}-\frac{1}{2}\rbr{m_i(1)-m_i(2)}}^2}d\bar{u}d\bar{v}\\
&=\frac{1}{2\pi}\frac{d}{\sigma^2}\frac{1}{2\sqrt{k}\sqrt{n-k}}\int_{\abs{u}\ls\va}\int_{\abs{v}\ls\va}
\exp\rbr{-\frac{d}{8\sigma^2 k}\abs{(u-m_i(1))+(v-m_i(2))}^2}\\
&\cdot \exp\rbr{-\frac{d}{8\sigma^2 (n-k)}\abs{(u-m_i(1))-(v-m_i(2))}^2}dudv\\
&=\frac{1}{2\pi}\frac{d}{\sigma^2}\frac{1}{2\sqrt{k}\sqrt{n-k}}\int_{\abs{u}\ls\va}\int_{\abs{v}\ls\va}
\exp\rbr{-\frac{dn}{8\sigma^2 k(n-k)}\rbr{\abs{u-m_i(1)}^2+\abs{v-m_i(2)}^2}}\\
&\cdot\exp\rbr{-\frac{d(n-2k)}{4\sigma^2 nk(n-k)}(u-m_i(1))(v-m_i(2))}dudv.
\end{align*}
The main problem in the above formula is $m_i(1)m_i(2)$. Before, we show a suitable upper bound, note that
\begin{align*}
&\frac{dn}{8\sigma^2 k(n-k)}\rbr{\abs{u-m_i(1)}^2+\abs{v-m_i(2)}^2}=\frac{d}{2\sigma^2 n}\rbr{\abs{u-m_i(1)}^2+\abs{v-m_i(2)}^2}\\
&-\frac{d\rbr{n-2k}^2}{8\sigma^2 nk(n-k)}\rbr{\abs{u-m_i(1)}^2+\abs{v-m_i(2)}^2}.
\end{align*}
Moreover,
\[
\frac{d\rbr{n-2k}^2}{8\sigma^2 nk(n-k)}\rbr{\abs{u-m_i(1)}^2+\abs{v-m_i(2)}^2}\gs \frac{d(n-2k)^2}{4\sigma^2 nk(n-k)} \abs{u-m_i(1)}\abs{v-m_i(2)}.
\]
Note that of $k\ls n/2$,
\[
\frac{d\abs{n-2k}n}{4\sigma^2 nk(n-k)}- \frac{d(n-2k)^2}{4\sigma^2 nk(n-k)} =\frac{d\abs{n-2k}}{2\sigma^2 n(n-k) }= \frac{d\abs{n-2k}}{2\sigma^2 n\max\cbr{k,n-k}}.
\]
Similarly, for $k>n/2$
\[
\frac{d\abs{n-2k} n}{4\sigma^2 nk(n-k)}- \frac{d(n-2k)^2}{4\sigma^2 nk(n-k)} =\frac{d\abs{n-2k}}{2\sigma^2 nk }= \frac{d\abs{n-2k}}{2\sigma^2 n\max\cbr{k,n-k}}.
\]
Therefore
\begin{align*}
& \frac{d(2k-n)}{4\sigma^2 k(n-k)}(u-m_i(1))(v-m_i(2))=\mathrm{sgn}\rbr{2k-n} \frac{d\abs{n-2k}}{\sigma^2 k(n-k)}(u-m_i(1))(v-m_i(2))\\
&=\mathrm{sgn}\rbr{2k-n}\rbr{\frac{d\abs{n-2k}}{2n\sigma^2\max\cbr{k,n-k}}+\frac{d(n-2k)^2}{4\sigma^2 nk(n-k)}}(u-m_i(1))(v-m_i(2))\\
&\ls \frac{d\rbr{2k-n}}{2n\sigma^2\max\cbr{k,n-k}}(u-m_i(1))(v-m_i(2))+\frac{d(n-2k)^2}{4\sigma^2 nk(n-k)}\abs{u-m_i(1)}\abs{v-m_i(2)}\\
& \ls \frac{d\rbr{2k-n}}{2\sigma^2 n\max\cbr{k,n-k}}(u-m_i(1))(v-m_i(2))+\frac{d(n-2k)^2}{8\sigma^2 n k(n-k)}\rbr{\abs{u-m_i(1)}^2+\abs{v-m_i(2)}^2}.
\end{align*}
It implies that
\begin{align*}
& \frac{1}{2\pi}\frac{d}{\sigma^2}\frac{1}{2\sqrt{k}\sqrt{n-k}}\int_{\abs{u}\ls\va}\int_{\abs{v}\ls\va}
\exp\rbr{ -\frac{dn}{8\sigma^2 k(n-k)}\rbr{\abs{u-m_i(1)}^2+\abs{v-m_i(2)}^2}}\\
&\cdot \exp\rbr{-\frac{d(n-2k)}{4\sigma^2 k(n-k)}(u-m_i(1))(v-m_i(2))}dudv \\
& \ls \frac{1}{2\pi}\frac{d}{\sigma^2}\frac{1}{2\sqrt{k}\sqrt{n-k}}\int_{\abs{u}\ls\va}\int_{\abs{v}\ls\va}
\exp\rbr{ \frac{d(n-2k)^2}{8\sigma^2 n k(n-k)}\rbr{\abs{u-m_i(1)}^2+\abs{v-m_i(2)}^2}}\\
& \cdot
\exp\rbr{-\frac{d}{2n \sigma^2}\rbr{\abs{u-m_i(1)}^2+\abs{v-m_i(2)}^2}}\\
& \cdot \exp\rbr{-\frac{d(n-2k)}{2\sigma^2 n\max\cbr{k,n-k}}(u-m_i(1))(v-m_i(2))}dudv.
\end{align*}
Clearly
\begin{align*}
& \exp\rbr{-\frac{d(n-2k)}{2\sigma^2 n\max\cbr{k,n-k}}(u-m_i(1))(v-m_i(2))}\\
& \ls \exp\rbr{-\frac{d(n-2k)}{2\sigma^2 n\max\cbr{k,n-k}}m_i(1)m_i(2)}\\
& \cdot \exp\rbr{\frac{d\abs{n-2k})}{2\sigma^2 n\max\cbr{k,n-k}}\rbr{\va\rbr{\abs{m_i(1)}+\abs{m_i(2)}}+\va^2}}.
\end{align*}
It shows that
\begin{align*}
& \frac{1}{2\pi}\frac{d}{\sigma^2}\frac{1}{2\sqrt{k}\sqrt{n-k}}\int_{\abs{u}\ls\va}\int_{\abs{v}\ls\va}
\exp\rbr{ -\frac{d}{2n \sigma^2}\rbr{\abs{u-m_i(1)}^2+\abs{v-m_i(2)}^2}}\\
& \cdot \exp\rbr{-\frac{d(n-2k)}{2\sigma^2 n\max\cbr{k,n-k}}(u-m_i(1))(v-m_i(2))}dudv\\
& \ls \frac{1}{2\pi}\frac{d}{\sigma^2}\frac{1}{2\sqrt{k}\sqrt{n-k}}
 \exp\rbr{\frac{d\abs{n-2k})}{2\sigma^2 n\max\cbr{k,n-k}}\rbr{\va\rbr{\abs{m_i(1)}+\abs{m_i(2)}}+\va^2}}\\
&\cdot \exp\rbr{-\frac{d(n-2k)}{2\sigma^2 n\max\cbr{k,n-k}}m_i(1)m_i(2)}\\
& \cdot \int_{\abs{u}\ls\va}\int_{\abs{v}\ls\va}\exp\rbr{ -\frac{d}{2\sigma^2 n}\rbr{\abs{u-m_i(1)}^2+\abs{v-m_i(2)}^2}}du dv
\end{align*}
Now we analyze the other side
\begin{align*}
&\P_{R^1}\rbr{\abs{\tbr{M_i+R^1_i,X}}\ls\va}\P_{R^2}\rbr{\abs{\tbr{M_i+R^2_i,Y}}\ls\va}\\
&=\P_{G^1_i}\rbr{\frac{\sigma}{\sqrt{d}}\abs{\sum^n_{j=1}G^1_{i,j}x_j-\sum^n_{j=1}M_{ij}x_j}\ls\va}\\
&\cdot\P_{G^2_i}\rbr{\abs{\frac{\sigma}{\sqrt{d}}\sum^n_{j=1}G^2_{ij}y_j-\sum^n_{j=1}m_{ij}y_j}\ls\va}\\
&=\P_{g_1}\rbr{\abs{\frac{\sigma}{\sqrt{d}}\sqrt{n}g_1-m_i(1))}\ls\va}\P_{g_2}\rbr{\abs{\frac{\sigma}{\sqrt{d}}\sqrt{n}g_2-m_i(2)}\ls\va}\\
&=\frac{1}{2\pi}\frac{d}{n\sigma^2}\int_{\abs{u}\ls\va}\int_{\abs{v}\ls\va} \exp\rbr{-\frac{d}{2 \sigma^2 n}\rbr{\abs{u-m_i(1)}^2+
\abs{v-m_i(2)}^2} }dudv.
\end{align*}
Therefore
\begin{align*}
& C_i(x,y)=\frac{n}{2\sqrt{k(n-k)}}\exp\rbr{\frac{d\abs{n-2k})}{2n\sigma^2\max\cbr{k,n-k}}\rbr{\va\rbr{\abs{m_i(1)}+\abs{m_i(2)}}+\va^2}}\\
&\cdot \exp\rbr{-\frac{d(n-2k)}{2n\sigma^2\max\cbr{k,n-k}}m_i(1)m_i(2)}.
\end{align*}
 Thus finally
\begin{align*}
& C(x,y)=\prod^d_{i=1}C_i(x,y)\ls \rbr{\frac{n}{2\sqrt{k(n-k)}}}^d \exp\rbr{-\frac{d(n-2k)}{2n\sigma^2\max\cbr{k,n-k}}\sum^d_{i=1}m_i(1)m_i(2)} \\
& \cdot 
\exp\rbr{\frac{d\abs{n-2k}}{2n\max\cbr{k,n-k}\sigma^2}\va \sum^d_{i=1}\rbr{\abs{m_i(1)}+\abs{m_i(2)}+\va}}.
 \end{align*}
Recall that for $x,y\in A$
\[
\sum^d_{i=1} m^2_i(1)\ls 2CV,\;\; \sum^d_{i=1} m^2_i(2)\ls 2CV.
\]
and thus
			 
\[
\sum^d_{i=1}\rbr{m_i(1)+m_i(2)}\ls \sqrt{d}\sqrt{8C V}.
\]
To go further we have to analyze $k=\abs{x+y}^2/4$, $n-k=\abs{x-y}^2/4$, $\tbr{x,y}=2k-n$.
In particular,
\[
k(n-k)=4^{-1}(n+\tbr{x,y})(n-\tbr{x,y})=4^{-1}\rbr{n^2-\tbr{x,y}^2},
\]
so
\[
\rbr{\frac{n}{2\sqrt{k(n-k)}}}^d=\rbr{\frac{1}{\sqrt{1-n^{-2}\tbr{x,y}^2}}}^d.
\]
Moreover, $2k-n=\tbr{x,y}$ and $\max\cbr{k,n-k}\gs n/2$ and
\[
\sum^d_{i=1}m_i(1)m_i(2)=\sum^d_{i=1}\tbr{M x,e_i}\tbr{M y,e_i}=\tbr{\bar{x},y},
\]
where $\hat{x}=\sum^d_{i=1}\tbr{Mx,e_i}e_i$ and observe that for $\norm{Mx}_2^2\ls 2 CV$ we have
\begin{equation}\label{ineq:3.3}
\norm{\hat{x}}_2^2=\sum^d_{i=1}\tbr{M,xe_i}^2\ls 2 CV.
\end{equation}
Hence
\begin{equation}\label{ineq3.3}
 C(x,y)\ls\rbr{\frac{1}{\sqrt{1-n^{-2}\tbr{x,y}^2}}}^d
 \exp\rbr{\frac{d\abs{\tbr{x,y}}}{n^2\sigma^2}\rbr{\va\sqrt{d}\sqrt{8 CV}+\va^2 d+\abs{\tbr{\bar{x},y}}}}.
\end{equation}

\section{How to split the bound}

Note that
\begin{align*}
&\P_{G_i}\rbr{\abs{\sum^n_{j=1}\frac{\sigma}{\sqrt{d}}G_{ij}x_j-\sum^n_{j=1}m_{ij}x_j}\ls\va}=\P_{g}\rbr{\abs{\frac{\sigma}{\sqrt{d}}\sqrt{n}g-\tbr{M_i,x}}}\\
&=\frac{\sqrt{d}}{\sqrt{2\pi n}\sigma}\int^{\va}_{-\va}
\exp\rbr{-\frac{d\rbr{u-\tbr{M_i,x}}^2}{2\sigma^2 n}}du\gs\exp\rbr{-\frac{d\tbr{M_i,x}^2}{2\sigma^2 n}}\gamma_1(\frac{\sqrt{d}\va}{\sqrt{n}\sigma}K_1),
\end{align*}
where $K_1=[-1,1]$ and $\gamma_1$ - standard Gaussian measure. Therefore, for $x\in A$
\begin{align*}
& \P_R\rbr{\norm{\rbr{M+R}x}_{\infty}\ls\va}=\prod^d_{i=1}\P_{g_i}\rbr{\abs{\frac{\sigma}{\sqrt{d}}\sqrt{n}g_i-\tbr{M_i,x}}\ls\va}\\
& \gs \exp\rbr{-\sum^d_{i=1}\frac{d\tbr{M_i,x}^2}{2\sigma^2 n}}\gamma_d(\frac{\sqrt{d}\va}{\sqrt{n}\sigma}K_d)
=\exp\rbr{-\frac{d\norm{Mx}_2^2}{2\sigma^2 n}}\gamma_d(\frac{\sqrt{d} \va}{\sqrt{n}\sigma}K_d),
\end{align*}
where $K_d=[-1,1]^d$ and $\gamma_d$ is the standard Gaussian measure on $\R^d$. Thus
\begin{align}\label{eq:4.1}
\begin{split}
& \E_Y \P_R\rbr{\norm{\rbr{M+R}Y}_{\infty}\ls\va}=W^{-1}\E_{\bar{Y}}\exp\rbr{\frac{d\norm{M\bar{Y}}_2^2}{2\sigma^2 n}}\1_{Y\in B}
\P_R\rbr{\norm{\rbr{M+R}\bar{Y}}_{\infty}\ls\va}\\
& \gs W^{-1}\P_{\bar{Y}}\rbr{Y\in B} \gamma_d(\frac{\sqrt{d} \va}{\sqrt{n}\sigma}K_d)\gs W^{-1} \rbr{1-C^{-1}}\gamma_d\rbr{\frac{\sqrt{d} \va}{\sqrt{n}\sigma}K_d},
\end{split}
\end{align}
where we used that $\P_{\bar{Y}}(\bar{Y}\in B)\gs 1-C^{-1}$ and that
\[
W=\E_{\bar{Y}} \1_{\bar{Y}\in B}\exp\rbr{\frac{d\norm{M\bar{Y}}_2^2}{2\sigma^2 n}}.
\]
On the other hand, 
\begin{align*}
& \frac{\sqrt{d}}{\sqrt{2\pi n}\sigma}\int^{\va}_{-\va}
\exp\rbr{-\frac{d\rbr{u-\tbr{M_i,x}}^2}{2\sigma^2 n}}du\\
&= \exp\rbr{-d\frac{\tbr{M_i,x}^2}{2n\sigma^2}} \frac{\sqrt{d}}{\sqrt{2\pi n}\sigma}\int^{\va}_{-\va}
\exp\rbr{-\frac{d\rbr{u}^2}{2n\sigma^2}}\cosh\rbr{\frac{du\tbr{M_i,x}}{\sigma^2 n}}du\\
& \ls \exp\rbr{-d\frac{\tbr{M_i,x}^2}{2\sigma^2 n}} \gamma_1\rbr{\frac{\sqrt{d} \va}{\sqrt{n}\sigma}}K_1)\cosh\rbr{\frac{d\va\tbr{M_i,x}}{\sigma^2 n}}.
\end{align*}
Consequently, we have on $\bar{Y}\in B$, i.e. $\norm{M\bar{Y}}\ls 2CV$,
\begin{align*}
&\P_R\rbr{\norm{\rbr{M+R}x}_{\infty}\ls\va}=\prod^d_{i=1}\P_{g_i}\rbr{\abs{\frac{\sigma}{\sqrt{d}}\sqrt{n}g_i-\tbr{M_i,x}}\ls\va}\\
&\ls\exp\rbr{-\sum^d_{i=1}\frac{d\tbr{M_i,x}^2}{2n\sigma^2}}\gamma_d\rbr{\frac{\sqrt{n}\sigma\va}{\sqrt{d}}K_d}\prod^d_{i=1}\cosh\rbr{\frac{d\va\tbr{M_i,x}}{n\sigma^2}}\\
& \ls \exp\rbr{ \rbr{\frac{d\va^2}{\sigma^2 n} -1}\frac{d\norm{Mx}_2^2}{2\sigma^2 n}} \gamma_d\rbr{\frac{\sqrt{d} \va}{\sqrt{n}\sigma}K_d}\ls \exp\rbr{\frac{d^2\va^2 C V}{\sigma^4 n^2}}\exp\rbr{-\frac{d\norm{MY}_2^2}{2\sigma^2 n}}\gamma_d\rbr{\frac{\sqrt{d} \va}{\sqrt{n}\sigma}K_d}.
\end{align*}
It shows that
\begin{align}
& \E_Y \P_R\rbr{\norm{\rbr{M+R}Y}_{\infty}\ls\va}= W^{-1}\E_{\bar{Y}} \exp\rbr{\frac{d\norm{M\bar{Y}}_2^2}{2\sigma^2 n}}\1_{\norm{M\bar{Y}}_2^2\ls 2 C V} 
\P_R\rbr{\norm{\rbr{M+R}\bar{Y}}_{\infty}\ls\va}\nonumber\\
& \ls W^{-1}\E_{\bar{Y}} \1_{Y\in B} \gamma_d\rbr{\frac{\sqrt{d} \va}{\sqrt{n}\sigma}K_d} \ls W^{-1} \exp\rbr{\frac{d^2\va^2 2 C V}{2\sigma^4 n^2}}\gamma_d\rbr{\frac{\sqrt{d} \va}{\sqrt{n}\sigma}K_d}.\label{eq:4.3}
\end{align}
Now we split the main bound
\begin{align*}
& \E_{R}\sbr{\P_X\rbr{\norm{\rbr{M+R}X}_{\infty}\ls\va}}^2\\
& =\E_{X,Y}\P_R\rbr{\cbr{\norm{\rbr{M+R}X}_{\infty}\ls\va}\cap\cbr{\norm{\rbr{M+R}Y}_{\infty}\ls\va}}\\
& =\E_{X,Y}\1_{\abs{\tbr{X,Y}}\ls n/2}\P_R\rbr{\cbr{\norm{\rbr{M+R}X}_{\infty}\ls\va}\cap\cbr{\norm{\rbr{M+R}Y}_{\infty}\ls\va}}\\
& + \E_{X,Y}\1_{\abs{\tbr{X,Y}}>n/2}\P_R\rbr{\cbr{\norm{\rbr{M+R}X}_{\infty}\ls\va}\cap\cbr{\norm{\rbr{M+R}Y}_{\infty}\ls\va}}.
\end{align*}

\section{Bound for $\abs{\tbr{X,Y}} > n/2$}

We first bound the case of $\abs{\tbr{x,y}} \gs n/2$. Note that
$\norm{M\bar{Y}}^2_2\ls 2CV$ for $\bar{Y}\in B$, so
\begin{align*}
\E_Y \1_{\abs{\tbr{X,Y}}>n/2}&=\E_{\bar{Y}}\1_{\abs{\tbr{X,\bar{Y}}}>n/2}\1_{\bar{Y}\in B}\exp\rbr{\frac{d\norm{M\bar{Y}}_2^2}{2\sigma^2 n}}\\
&\ls\exp\rbr{\frac{d C V}{\sigma^2 n}} \P_{\bar{Y}} \rbr{\abs{\tbr{X,\bar{Y}}}>n/2}.
\end{align*}
This event has a small probability, namely
\begin{align*}
\P_{\bar{Y}}\rbr{\abs{\tbr{X,\bar{Y}}}>n/2}&\ls 2e^{-\lambda n/2}\E_{\bar{Y}}\exp\rbr{\lambda\tbr{X,\bar{Y}}}=2e^{-\frac{1}{2}\lambda n+\lambda^2 n}\\
&\ls 2e^{-n/16},\;\;\mbox{for}\;\;\lambda=1/4.
\end{align*}
Therefore, by \eqref{eq:4.1} we have
\begin{align*}
& \E_{X,Y} \1_{\abs{\tbr{X,Y}}>n/2}\P_R\rbr{\cbr{\norm{\rbr{M+R}X}_{\infty}\ls\va}]\cap \cbr{\rbr{\norm{\rbr{M+R}Y}_{\infty}\ls\va}}} \\
& \ls \E_{X,Y} \1_{\abs{\tbr{X,Y}}>n/2}\P_R\rbr{\norm{\rbr{M+R}X}_{\infty}\ls\va} \\
&= \E_X \P_R\rbr{\norm{\rbr{M+R}X}_{\infty}\ls\va}\E_Y \1_{\abs{\tbr{X,Y}}>n/2}\\
& \ls 2e^{-n/16}\exp\rbr{\frac{d C V}{\sigma^2 n}} \E_X \P_R\rbr{\norm{\rbr{M+R}X}_{\infty}\ls\va}\\
& \ls 2e^{-n/16}\exp\rbr{\frac{d C V}{\sigma^2 n}} W \rbr{1-C^{-1}}^{-1} \gamma^{-1}_d
\rbr{\frac{\sqrt{d} \va}{\sqrt{n}\sigma}K_d}\sbr{ \E_X \P_R\rbr{\norm{\rbr{M+R}X}_{\infty}\ls\va}}^2.
\end{align*}
In order to bound $W$ we can use
\[
W\ls \exp\rbr{\frac{d C V}{\sigma^2 n}}, 
\]
and thus finally
\begin{align*}
& \E_{X,Y} \1_{\abs{\tbr{X,Y}}>n/2} \P_R\rbr{ \cbr{\norm{\rbr{M+R}X}_{\infty}\ls\va}]\cap \cbr{\rbr{\norm{\rbr{M+R}Y}_{\infty}\ls\va}}} \\
& \ls 2(1-C^{-1})^{-1}e^{-n/16}\exp\rbr{\frac{2d C V}{\sigma^2 n}} \gamma^{-1}_d
\rbr{\frac{\sqrt{d} \va}{\sqrt{n}\sigma}K_d}\sbr{ \E_X \P_R\rbr{\norm{\rbr{M+R}X}_{\infty}\ls\va}}^2.
\end{align*}
We aim to analyze when 
\[
 2e^{-n/16}\exp\rbr{\frac{2d CV}{\sigma^2 n}}\gamma^{-1}_d\rbr{\frac{\sqrt{d} \va}{\sqrt{n}\sigma}K_d}\ls 2e^{-n/64}.
\]
This can be done be requiring that
\begin{equation}\label{eq:4.4}
\gamma_d\rbr{\frac{\sqrt{d} \va}{\sqrt{n}\sigma}K_d}\gs e^{-n/32}
\end{equation}
and
\begin{equation}\label{eq:4.5}
\frac{2d CV}{\sigma^2 n}\ls \frac{n}{64},\;\; n\gs 8 \sqrt{d} \sqrt{2 C V}.
\end{equation}
Clearly \eqref{eq:4.4} reads as
\[
\gamma_1 \rbr{\frac{\sqrt{d} \va}{\sqrt{n}\sigma}K_1}\gs e^{-n/\rbr{32 d}},
\;\;\mbox{i.e.}\;\;\P\rbr{\abs{g}\ls \frac{ \sqrt{d}\va}{\sigma\sqrt{n} }}\gs e^{-\frac{n}{32 d}}.
\]
Therefore, we require that $\va\gs\va_0$ for $\va_0$ satisfying
\[
\P\rbr{\abs{g}\ls \frac{ \sqrt{d}\va_0}{\sigma\sqrt{n} }}= e^{-\frac{n}{32 d}}.
\]
We require that $\sqrt{d}\va_0/\rbr{\sigma \sqrt{n}}\ls c$ for a small positive constant $c$, and thus we can use the inequality
\[
P\rbr{\abs{g}\ls \frac{ \sqrt{d}\va_0}{\sigma\sqrt{n} }}\gs \frac{2\va_0\sqrt{d}}{\sqrt{2\pi}\sigma\sqrt{n}}\exp\rbr{-\frac{c^2}{2}}.
\]
Hence
\[
\va_0\ls \sqrt{\frac{\pi}{2}}e^{\frac{c^2}{2}}\frac{\sigma \sqrt{n}}{\sqrt{d}}e^{-\frac{n}{32d}}
\]
and thus we assume that
	 
\[
\va \gs \sqrt{\frac{\pi}{2}}e^{\frac{c^2}{2}}\frac{\sigma \sqrt{n}}{\sqrt{d}}e^{-\frac{n}{32d}}.
\]
The right interval for $\va$ is
\[
\frac{\sigma \sqrt{n}}{\sqrt{d}}e^{-\frac{n}{32 d}}\lesssim \va .
\]
Therefore, 
\begin{align*}
 & \E_{X,Y} \1_{\abs{\tbr{X,Y}}>n/2} \P_R\rbr{ \cbr{\norm{\rbr{M+R}X}_{\infty}\ls\va}]\cap \cbr{\rbr{\norm{\rbr{M+R}Y}_{\infty}\ls\va}}} \\
& \ls 2e^{-n/64} \rbr{1-C^{-1}}^{-1}\sbr{\E_X \P_R\rbr{\norm{\rbr{M+R}X}_{\infty}\ls\va}}^2.
\end{align*}

\section{Bound for $\tbr{X,Y} \ls n/2$}

Now we analyze the case when $\abs{\tbr{x,y}}<n/2$. Here we use the our bound on $C(x,y)$.
Namely, by \eqref{eq:4.3}, we have
\begin{align*}
& \E_{X,Y} \1_{\abs{\tbr{X,Y}}\ls n/2}\P_R\rbr{\cbr{\norm{\rbr{M+R}X}_{\infty}\ls\va}\cap \cbr{\norm{\rbr{M+R}Y}_{\infty}\ls\va}}\\
& \ls \E_{X,Y} \1_{\abs{\tbr{X,Y}}\ls n/2}C(X,Y) \P_R\rbr{\norm{\rbr{M+R}X}_{\infty}\ls\va}\P_R\rbr{\norm{\rbr{M+R}Y}_{\infty}\ls\va}\\
& \ls W^{-1}\exp\rbr{\frac{d^2\va^2 2 C V}{2\sigma^4 n^2}}\gamma_d\rbr{\frac{\sqrt{d}\va}{\sqrt{n}\sigma}K_d}. 
\E_X \P_{R} \rbr{\norm{\rbr{M+R}X}_{\infty}\ls\va}\E_{\bar{Y}}\1_{\abs{\tbr{X,\bar{Y}}}\ls n/2}C(X,\bar{Y}).
\end{align*}
Note that
\[
\frac{1}{\sqrt{1-x^2}}=\exp\rbr{-\frac{1}{2}\ln\rbr{1-x^2}}\ls e^{\frac{1}{2}x^2+\frac{1}{2}x^4}\ls \exp\rbr{\frac{5}{8}x^2},\;\;\mbox{for}\;\;\abs{x}<1/2.
\]
Recall the definition of $\hat{X}$ i.e. \eqref{ineq:3.3}. We have
\begin{align*}
& \E_{\bar{Y}} \1_{\abs{\tbr{X,\bar{Y}}}\ls n/2} C(X,\bar{Y})\\
&\ls \E_{\bar{Y}} \rbr{\frac{1}{\sqrt{1-n^{-2}\tbr{X,\bar{Y}}^2}}}^d
\exp\rbr{\frac{d\abs{\tbr{X,\bar{Y}}}}{\sigma^2 n^2}\rbr{\va^2 d+\va\sqrt{d}\sqrt{ 8CV} }+\abs{\tbr{\hat{X},\bar{Y}}}}\\
& \ls \E_{\bar{Y}} \exp\rbr{ \frac{5d}{8n^2}\tbr{X,\bar{Y}}^2 +\frac{d\abs{\tbr{X,\bar{Y}}})}{\sigma^2 n^2}
\rbr{\va^2 d+\va\sqrt{d}\sqrt{ 8CV} +\abs{\tbr{\hat{X},\bar{Y}}}}} .
\end{align*}
By Schwarz inequality
\begin{align*}
& \E_{\bar{Y}} \1_{\abs{\tbr{X,\bar{Y}}}\ls n/2} C(X,\bar{Y})\\
& \ls \rbr{\E_{\bar{Y}} \exp\rbr{ \frac{5d}{4n^2} \tbr{X,\bar{Y}}^2 +
\frac{2d\abs{\tbr{X,\bar{Y}}})}{\sigma^2 n^2}\rbr{\va^2 d+\va\sqrt{d}\sqrt{ 8CV}}}}^{1/2} \\
& \cdot \rbr{\E_{\bar{Y}} \exp\rbr{\frac{2d}{\sigma^2 n^2}\abs{\tbr{X,\bar{Y}}}\abs{\tbr{\hat{X},\bar{Y}}}}}^{1/2}. 
\end{align*}
However, $\abs{\tbr{X,\bar{Y}}}$ satisfies
\begin{align*}
\P_{\bar{Y}} \rbr{\abs{\tbr{X,\bar{Y}}} >u}\ls 2e^{-\lambda u}\E_{\bar Y} \exp\rbr{\lambda \tbr{X,\bar{Y}}}=2e^{-\lambda u+\lambda^2 n}\ls 2e^{-u^2/\rbr{4n}}
\end{align*}
for $\lambda=u/(2n)$.
Now for $c>b$
\begin{align*}
& \int^{\infty}_0 \rbr{2bt+a}e^{at+bt^2}e^{-ct^2}dt \\
& \int^{\infty}_0 \frac{-b}{c-b}\rbr{-2(c-b)t+a}e^{at-(c-b)t^2}dt+\frac{ca}{c-b} \int^{\infty}_0 e^{at}e^{-(c-b)t^2}dt\\
&= \frac{b}{c-b}+\frac{ca}{c-b}\frac{\sqrt{\pi}}{\sqrt{c-b}} \frac{\sqrt{2(c-b)}}{\sqrt{2\pi}}\int_{\R} e^{at}e^{-(2(c-b))t^2/2} dt\\
& =\frac{b}{c-b}+\frac{ca}{c-b}\frac{\sqrt{\pi}}{\sqrt{c-b}}e^{\frac{a^2}{4(c-b)}}.
\end{align*}
In our case
\[
a=\frac{2d}{\sigma^2 n^2}\rbr{\va^2 d+\va\sqrt{d}\sqrt{ 8CV}}, \;\; b= \frac{5d}{4n^2}\,\;\;c=\frac{1}{4n}
\]
The crucial assumption
\[
b\ls \frac{1}{2}c,\;\;\mbox{ie.}\;\; 10d\ls n
\]
and
\[
\frac{a}{\sqrt{b-c}}=\sqrt{8n}\frac{2d}{\sigma^2 n^2}\rbr{\va^2 d+\va\sqrt{d}\sqrt{ 8CV}}\ls \delta,
\]
where $\delta$ is suitably small.
It holds true, for example when 
\[
\va \ls \frac{\delta}{32}\sigma^2\rbr{\frac{n}{d}}^{3/2}(C V)^{-1/2},\;\;\va \ls \frac{\delta^{1/2}}{4}\sigma \rbr{\frac{n}{d}} n^{-1/4}.
\]
Therefore,
\[
 \frac{b}{c-b}+\frac{ca}{c-b}\frac{\sqrt{\pi}}{\sqrt{c-b}}e^{\frac{a^2}{4(c-b)}} \ls 10\frac{d}{n}+2\sqrt{\pi}\delta e^{\frac{\delta^2}{4}}
\]
It shows that
\begin{align*}
& \rbr{\E_{\bar{Y}} \exp\rbr{ \frac{d}{n^2}\rbr{\frac{5}{4}+\frac{1}{\sigma^2}} \tbr{X,\bar{Y}}^2 +
\frac{2d\abs{\tbr{X,\bar{Y}}})}{\sigma^2 n^2}\rbr{\va^2 d+\va\sqrt{d}\sqrt{ 8CV}}}}^{1/2} \\
& \ls \rbr{1+2\int^{\infty}_0 \rbr{2bt+a}e^{at+bt^2}e^{-ct^2}dt}^{1/2}\\
& \ls \rbr{1+ 20\frac{d}{n}+4\sqrt{\pi}\delta e^{\frac{\delta^2}{4}}}^{1/2}.
\end{align*}
On the other hand,
\begin{align*}
& \E_{\bar{Y}} \exp\rbr{\frac{2d}{\sigma^2 n^2}\abs{\tbr{X,\bar{Y}}}\abs{\tbr{\bar{X},\bar{Y}}}}\\
& \rbr{\E_{\bar{Y}} \exp\rbr{\frac{2d}{\sigma^2 n^2}\tbr{X,\bar{Y}}^2 }}^{1/2}
\rbr{\E_{\bar{Y}} \exp\rbr{\frac{2d}{\sigma^2 n^2}\tbr{\hat{X},\bar{Y}}^2 }}^{1/2}.
\end{align*}
We know that
\begin{align*}
&\rbr{\E_{\bar{Y}}\exp\rbr{\frac{2d}{\sigma^2 n^2}\tbr{X,\bar{Y}}^2 }}^{1/2}\\
&\ls \rbr{1+2\int^{\infty}_0 2bt e^{bt^2}e^{-ct^2}dt}^{1/2}=\rbr{1+2\frac{b}{c-b}}^{1/2},
\end{align*}
where
\[
b=\frac{2d}{\sigma^2 n^2},\;\;c=\frac{1}{4n}
\]
we need that
\[
b<\frac{1}{2}c,\;\; \mbox{i.e.}\;\; \frac{n}{d}\gs 16\sigma^2.
\]
Then
\[
 \rbr{\E_{\bar{Y}}\exp\rbr{\frac{2d}{\sigma^2 n^2}\tbr{X,\bar{Y}}^2}}^{1/2}\ls \rbr{1+32\frac{d}{n\sigma^2}}^{1/2}.
\]
Now observe that $\abs{\tbr{\hat{X},\bar{Y}}}$ satisfies
\begin{align*}
& \P_Y \rbr{\abs{\tbr{\hat{X},\bar{Y}}}>u}\ls 2e^{-\lambda u} \E_{\bar{Y}} \exp\rbr{\lambda \tbr{\hat{X},\bar{Y}} }\\
& =2e^{-\lambda u+\lambda^2 2C V}\ls 2e^{-u^2/\rbr{8 C V}},\;\;
\mbox{for}\;\;\lambda=\frac{u}{4 C V}.
\end{align*}
We can use our formula with
\[
b=\frac{2d}{\sigma^2 n^2},\;\;c=\frac{1}{8 C V}. 
\]
and
\[
b<c/2,\;\;n>\frac{4\sqrt{ d}\sqrt{C V}}{\sigma}.
\]
Then
\begin{align*}
\rbr{\E_{\bar{Y}}\exp\rbr{\frac{d}{\sigma^2 n^2}\tbr{\hat{X},\bar{Y}}^2}}^{1/2}&=\rbr{1+2\int^{\infty}_0 2bte^{bt^2}e^{-ct^2}dt }^{1/2}=\rbr{1+2\frac{b}{c-b}}^{1/2}\\
&\ls\rbr{1+\frac{32d CV}{\sigma^2n^2}}^{1/2}.
\end{align*}
It implies that
\begin{align*}
\rbr{\E_{\bar{Y}}\exp\rbr{\frac{2d}{\sigma^2 n^2}\abs{\tbr{X,\bar{Y}}}\abs{\tbr{\hat{X},\bar{Y}}}}}^{1/2}\ls\rbr{1+32\frac{d}{\sigma^2 n}}^{1/4}\rbr{1+ \frac{32d C V}{\sigma^2n^2}}^{1/4}.
\end{align*}
Hence
\begin{align}
\begin{split}
\E_{\bar{Y}}\1_{\abs{\tbr{X,\bar{Y}}}\ls n/2}C(X,\bar{Y})&\ls\rbr{1+ 20\frac{d}{n}+4\sqrt{\pi}\delta e^{\frac{\delta^2}{4}}}^{1/2} \rbr{1+32\frac{d}{n\sigma^2}}^{1/4}\rbr{1+ \frac{32d C V}{\sigma^2n^2}}^{1/4}\\
&=1+\bar{\delta}.
\end{split}
\end{align}
We collect conditions under which we may claim that $\bar{\delta}$ is suitably small. We need that
\begin{align}
& \frac{n}{d}> L,\;\;\frac{n}{d}>\frac{L}{\sigma^2};\nonumber\\
& \va \ls \frac{\delta}{32}\sigma^2\rbr{\frac{n}{d}}^{3/2}(C V)^{-1/2},\;\;\va \ls \frac{\delta^{1/2}}{4}\sigma \rbr{\frac{n}{d}} n^{-1/4},\label{eq:4.6}
\end{align}
where $L$ is suitably large and $\delta$ is suitably small.
Thus by \eqref{eq:4.3}
\begin{align*}
& \E_{X,Y} \1_{\abs{\tbr{X,Y}}\ls n/2}
\P_R\rbr{\cbr{\norm{\rbr{M+R}X}_{\infty}\ls\va}\cap \cbr{\norm{\rbr{M+R}Y}_{\infty}\ls\va}}\\
& \ls \rbr{1+\bar{\delta}}\exp\rbr{\frac{d^2\va^2 2 C V}{2\sigma^4 n^2}} W^{-1} \gamma_d(\frac{\sqrt{n}\sigma \va}{\sqrt{d}}K_d) 
\E_X \1_{X\in A} \P_{R^1}\rbr{\norm{\rbr{M+R^1}X}_{\infty}\ls\va}\\
&\ls \rbr{1+\bar{\delta}}\rbr{1-C^{-1}}^{-1}\exp\rbr{\frac{d^2\va^2 C V}{\sigma^4 n^2}}\sbr{\E_X \1_{X\in A} \P_R\rbr{\norm{\rbr{M+R}X}_{\infty}\ls\va}}^2.
\end{align*}
Note that by our assumptions
\[
\va\ls\frac{\delta}{32}\sigma^2\frac{n}{d}(CV)^{-1/2},\;\;\mbox{so}\;\;\exp\rbr{\frac{d^2\va^2 CV}{\sigma^4 n^2}}\ls\exp\rbr{\frac{\delta^2}{32}}.
\]
Thus
\begin{align*}
& \E_{X,Y} \1_{\abs{\tbr{X,Y}}\ls n/2}\1_{X\in A}\1_{Y\in A} \P_R\rbr{\norm{\rbr{M+R}X}_{\infty}\ls\va}\P_{R}\rbr{\norm{\rbr{M+R}Y}_{\infty}\ls\va}\\
& \ls \rbr{1+\bar{\delta}}\gamma_d(\frac{\sqrt{n}\sigma \va}{\sqrt{d}}K_d)
\E_X \1_{X\in A} \P_{R^1} \rbr{\norm{\rbr{M+R^1}X}_{\infty}\ls\va}\\
&\ls (1+\bar{\delta})(1-C^{-1})^{-1}\sbr{\E_X \1_{X\in A} \P_R\rbr{\norm{\rbr{M+R}X}_{\infty}\ls\va}}^2.
\end{align*}
Thus matching our bounds
\begin{align*}
& \E_{R} \sbr{\P_X\rbr{\cbr{\norm{\rbr{M+R}X}_{\infty}\ls\va}\cap \cbr{X\in A}}}^2\\
&\ls \rbr{1+\bar{\delta}+e^{-n/32}}(1-C^{-1})^{-1}\sbr{\E_X \1_{X\in A} \P_R\rbr{\norm{\rbr{M+R}X}_{\infty}\ls\va}}^2.
\end{align*}
So
\[
\P_R \rbr{P_X\rbr{\cbr{\norm{\rbr{M+R}X}_{\infty}\ls\va}\cap \cbr{X\in A}}>0}\gs \frac{1}{1+\bar{\delta}+2e^{-n/32}}(1-C^{-1})
\]
Finally, we have to compare our requirements for $\va$. We consider $n=\kappa d\log d$.
We note that
\[
\va\gs\sqrt{\frac{\pi}{2}}e^{\frac{c^2}{2}}\frac{\sigma\sqrt{n}}{\sqrt{d}}e^{-\frac{n}{32d}}\sim\sigma\sqrt{\log d}d^{-\kappa/32}.
\]
Moreover, for $C\ls\rbr{\log{d}}^2$ and $\sigma\gs 1$, due to \eqref{property3}, the following inequalities hold true for $\delta\gs 32d^{1/2-\kappa/32}$
\[
\frac{\delta}{32}\sigma^2\rbr{\frac{n}{d}}^{3/2}(CV)^{-1/2}=\frac{\delta}{32}\sigma^2\rbr{\log{d}}^{3/2}(CV)^{-1/2}\gs\sigma\sqrt{\log{d}}d^{-\kappa/32}
\]
\[
\frac{\delta^{1/2}}{4}\sigma\rbr{\frac{n}{d}}n^{-1/4}=\frac{\delta^{1/2}}{4}\sigma\rbr{\log{d}}^{3/4}d^{-1/4}\gs\sqrt{2}\sigma\rbr{\log{d}}^{3/4}d^{-\kappa/64}\gs\sigma\sqrt{\log{d}}d^{-\kappa/32}.
\]
Thus conditions for $\varepsilon$ in \eqref{eq:4.6} are satisfied. This shows that Theorem \ref{theo2} holds true for $n=\kappa d\log d$, with $\va=\sigma \sqrt{\log d}d^{-\kappa/32}$.  It is also worth to note that 
additionally assuming that $C=\rbr{\log{d}}^2$,  $\delta=32 d^{1/2-\kappa/32}$ one can show that
$c(d)\sim 1/\log d$.


\begin{thebibliography}{99}

\footnotesize

\bibitem{Ban} \textsc{Banaszczyk, W.} (1998) Balancing vectors and gaussian measures of n-dimensional convex
bodies. \textit{Random Structures and Algorithms}, {\bf 12}(4):351–360.
\bibitem{Bansal} \textsc{Bansal, N., Dadush, D., Garg, S. and Lovett, S.} (2019). The Gram-SchmidtWalk: A cure for the Banaszczyk
Blues. \textit{Theory of Computing} ,{\bf 15}, (21), 1–27.
\bibitem{Nik} \textsc{Nikolov, A.} (2013) The Komlos conjecture holds for vector colorings. Available online at
\textit{arXiv:1301.4039}.
\bibitem{Meka} \textsc{Bansal, N., Jiang, H., Meka, R., Singla, S. and Sinha, M.} (2022) Smoothed Analysis of the Komlós Conjecture
\textit{arXiv:2204.11427}.
\bibitem{Spe1} \textsc{Spencer, J.} (1985) Six standard deviations suffice. Trans. Amer. Math. Soc., {\bf 289}, 679–679..
\bibitem{Spe2} \textsc{Spencer, J.} (1994) Ten Lectures on the Probabilistic Method: Second Edition. \textit{SIAM}.
\bibitem{Spiel} \textsc{Harshaw, C., Sävje, F., Spielman, D.A. and Zhang, P.} (2022) Balancing covariates in randomized experiments
with the Gram–Schmidt Walk design, \textit{arXiv:1911.03071}.

\end{thebibliography}
\end{document}